\documentclass[12 pt]{article}
\usepackage{amsfonts}
\title{Some Properties of Alphabet Overlap Graphs}
\author{Anant P.~Godbole, Debra Knisley, and Rick Norwood\\
Department of Mathematics\\
East Tennessee State University}
\begin{document}
\def\qed{\vbox{\hrule\hbox{\vrule\kern3pt\vbox{\kern6pt}\kern3pt\vrule}\hrule}}
\def\ep{\varepsilon}
\def\lr{\left(}
\def\lf{\lfloor}
\def\rf{\rfloor}
\def\lc{\left\{}
\def\rc{\right\}}
\def\rr{\right)}
\def\p{\mathbb P}
\def\v{\mathbb V}
\def\e{\mathbb E}
\def\l{\mathbb L}
\def\lg{{\rm lg}}
\maketitle
\begin{abstract}
Consider a graph $G=G(k,d,s)$ with the vertex set $V=\{v: v=(v_1,\ldots,v_k); v_i\in\{1,2,\ldots, d\}(1\le i\le k)\}$, the set of all $k$-letter ``words" over an ``alphabet" of size $d$.  Furthermore, there will be an edge between vertices $v\ne w$ iff the last $k-s$ letters of $v$ are the same as the first $k-s$ letters of $w$ {\it or} the first $k-s$ letters of $v$ are the same as the last $k-s$ letters of $w$.  In this paper, we show that $G$ is Hamiltonian for all non-trivial values of the parameters, and obtain exact values for its chromatic number when $s\ge k/2$ and bounds on its chromatic number when $s<k/2$.
\end{abstract}

\nocite{godbole}
\nocite{dna}
\nocite{collins}
\nocite{clark}
\nocite{diestel}
\nocite{robin}
\nocite{cameron}
\nocite{shields}
\def\qed{\vbox{\hrule\hbox{\vrule\kern3pt\vbox{\kern6pt}\kern3pt\vrule}\hrule}}
\def\ep{\varepsilon}
\def\lr{\left(}
\def\lf{\lfloor}
\def\rf{\rfloor}
\def\lc{\left\{}
\def\rc{\right\}}
\def\rr{\right)}
\def\p{\mathbb P}
\def\v{\mathbb V}
\def\e{\mathbb E}
\def\l{\mathbb L}
\def\lg{{\rm lg}}
\newtheorem{thm}{Theorem}
\newtheorem{lemma}[thm]{Lemma}
\newtheorem{prop}[thm]{Proposition}
\maketitle

\begin{abstract}
Consider a graph $G=G(k,d,s)$ with the vertex set $V=\{v: v=(v_1,\ldots,v_k); v_i\in\{1,2,\ldots, d\}(1\le i\le k)\}$, the set of all $k$-letter ``words" over an ``alphabet" of size $d$.  Furthermore, there will be an edge between vertices $v\ne w$ iff the last $k-s$ letters of $v$ are the same as the first $k-s$ letters of $w$ {\it or} the first $k-s$ letters of $v$ are the same as the last $k-s$ letters of $w$.  In this paper, we show that $G$ is Hamiltonian for all non-trivial values of the parameters, and obtain exact values for its chromatic number when $s\ge k/2$ and bounds on its chromatic number when $s<k/2$.
\end{abstract}

\nocite{godbole}
\nocite{dna}
\nocite{collins}
\nocite{clark}
\nocite{diestel}
\nocite{robin}
\nocite{cameron}
\nocite{shields}
\section{Introduction}
Consider a graph $G=G(k,d,s)$ with the vertex set
$$V=\{v: v=(v_1,\ldots,v_k); v_i\in\{1,2,\ldots, d\}(1\le i\le k)\},$$
the set of all $k$-letter ``words" over an ``alphabet" of size $d$.  We shall refer to the coordinates of $v$ as ``letters," and draw the reader's attention to the cases $d=2, d=4, d=20,$ and $d=26$ as providing concrete applications to binary words; DNA sequences; protein sequences; and words from the English language respectively.  Furthermore, there is an edge between vertices $v$ and $w$ iff $v\ne w$ and the last $k-s$ letters of $v$ are the same as the first $k-s$ letters of $w$ {\it or} the first $k-s$ letters of $v$ are the same as the last $k-s$ letters of $w$.  In this paper, we

(a) show that $G$ is Hamiltonian for all non-trivial values of the parameters $k,d,$ and $s$, and

(b) find exact values for the chromatic number $\chi(G)$ when $s\ge k/2$ and bounds on $\chi(G)$ when $s< k/2$.

\noindent (The result in part (a) is well known when $s=1$ and can easily be extended to the case where $s\le k/2$, so the novelty lies in the method of proof and the case where $s>k/2$.)

  We exhibit similarities and distinctions between $G$ and
(i) the standard $k$-cube and (ii) the $k$-dimensional grid on
$d^k$ points.  We call $G$ an {\it alphabet overlap graph}, noting
that the authors of \cite{dna} have studied similar graphs --
calling them $(\alpha,k)$ labeled graphs (their ``$\alpha$" is the
same as our ``$d$").  In their nomenclature, our graphs would be
best termed {\it complete $(d,k)$ labeled graphs}.  The motivation
in \cite{dna} was (i) to investigate connections between families
of such graphs with different parameter values; (ii) to develop
recognition algorithms; and (iii) to consider the case
$d\to\infty$, whereas we are more interested in the structural
properties of alphabet overlap graphs.

Discrete mathematics, graph theory in particular, is playing an
ever increasing role in the science of molecular biology. This is
evidenced, in part, by DIMACS's ``Special Years of Focus on
Computational Biology" (2000-2003); see \cite{dimacs} for details.
With the increased quantity and complexity of biodata, new tools
and frameworks are being developed to recognize and understand the
many processes involved. Graph theoretical approaches which employ
such tools as minimum spanning trees and bipartite matchings are
appearing more frequently in the literature. The whole field
representing the interplay between graph theory and molecular
biology thus appears to be one of the more exciting emerging areas
of interdisciplinary research, as further evidenced by, e.g., the
Special Session on Applications of Graph Theory in Molecular
Biology held at the 2002 SIAM Discrete Mathematics Meeting in San
Diego, and chaired by the second-named author of this paper.
Earlier work on the interface between graph theory and biology may
be found in \cite{roberts}.

Researchers in discrete probability have long studied {\it random letter generation} from a $d$-letter alphabet; see, e.g., \cite{bala} or the collection of papers in \cite{godbole}.  In this case, a Markov chain model is often appropriate with a ``transition" from the state $v=(v_1,\ldots,v_k)$ to $w=(v_2,\ldots,v_k,w_1)$ occuring, for {\it each} value of $w_1$, with probability $1/d$.  (In our context, however, the stated caveat that $v\ne w$ does somewhat more than simply reflect that fact that our graph has no loops; it states that  probabilistically feasible transitions such as those from AAAAA to AAAAA will not be allowed in the graph structure). If we randomly generate $n$ letters from a $d$-letter alphabet, then there are $n-k$ transitions of the type described above, and researchers have studied random variables such as $X$, defined as the number of occurrences, with overlaps possibly allowed, of fixed words such as ``ABRACADABRA".  The literature is replete with results along these lines, complete with distributional approximations, multivariate analogs, and connections with the overlap structure of the word in question.  See \cite{shields} for results of a different kind and \cite{guibas} for a landmark paper than provides necessary and sufficient conditions for a set $A$ of integers to be the set of periods of some string.

From a lay point of view, our main result on Hamiltonicity alluded to above states, for example, that we may start with any eight letter English word such as CATACOMB, make a ``transition" to ATACOMBO, and continue to make transitions so that no word is obtained twice; each of the $26^8$ words of length 8 are covered; and the starting word CATACOMB is recovered at the end of the Hamiltonian cycle.  The same is true of any longer starting document such as President George W. Bush's inaugural address, or, using the standard keyboard as an alphabet, even the {\tt .tex} file corresponding to a long mathematical paper.  From a technical point of view, we note that when $s=1$, the Hamiltonicity of $G$ is equivalent to the existence of a deBruijn cycle on the set of all $k$-letter words on a $d$-letter alphabet, which in turn may be viewed as an efficient ordering of this set.  An entire workshop was held recently at the Banff International Research Station that focused exclusively on deBruijn cycles and Gray codes; see \cite{birs} for details.

The standard sufficient conditions for Hamiltonicity due to Dirac et al.(see, e.g., Chapter 10 in \cite{diestel}) are simply not valid in our context, nor are more recent criteria such as the one in \cite{robin}. Of particular note is the fact that our proof is elementary, and uses induction on the alphabet size and not the word length (as is done, e.g., in \cite{cameron} for the $k$-cube.)

There are important differences between the geometry or
``architecture" of the $k$-cube $Q^k$ and the $k$-grid
$\{1,2,\ldots,d\}^k$ on the one hand, and that of alphabet overlap
graphs on the other.  In the former case, there is an edge between
vertices $v$ and $w$ iff $h(v,w)=1$, where $h$ represents Hamming
distance.  On the other hand, for small values of the parameters
our graphs are more-or-less like ``twisted" cubes or grids; see
Figure 1 for a drawing of $G=G(k,d,s)=G(3,2,1)$ from which a
comparison with $Q^3$ can be readily made.  The adjective
``twisted" should be used with caution, however, when we look at
larger structures.  If, for example, $s=1$, it is clear that ${\rm
deg}(v_i)\le 2d$ for each $i$, and that the inequality could be
strict:  the degree of words such as AAA...A is $2d-2$, whereas
that of words such as ABAB...$\diamondsuit$ is $2d-1$, where
$\diamondsuit$ is either A or B.  It follows that the number of
edges in alphabet overlap graphs with $s=1$ equals
\begin{eqnarray*}&&{}\lr{{2d(d^k-d-d(d-1))+(2d-1)(d(d-1))+(2d-2)d}\over{2}}\rr\nonumber\\
&&{}=d^{k+1}+O(d^3)\enspace(d\to\infty).
\end{eqnarray*}

\medskip

\setlength{\unitlength}{.2in}
\begin{figure}{}
\begin{center}
\begin{picture}(6,8)
\put(0,1){\circle*{.20}} \put(0,4){\circle*{.20}}
\put(3,1){\circle*{.20}} \put(3,4){\circle*{.20}}
\put(2,5){\circle*{.20}} \put(5,5){\circle*{.20}}
\put(5,8){\circle*{.20}} \put(2,8){\circle*{.20}}
\put(8,1){\circle*{.20}} \put(11,1){\circle*{.20}}
\put(8,4){\circle*{.20}} \put(11,4){\circle*{.20}}
\put(10,5){\circle*{.20}}

\put(13,5){\circle*{.20}}

\put(13,8){\circle*{.20}}

\put(10,8){\circle*{.20}}

\put(0,1){\line(3,0){3}}

\put(0,1){\line(0,3){3}}

\put(0,4){\line(3,0){3}}

\put(3,1){\line(0,3){3}}

\put(3,1){\line(1,2){2}}

\put(0,1){\line(1,2){2}}

\put(0,4){\line(1,2){2}}

\put(3,4){\line(1,2){2}}

\put(2,5){\line(3,0){3}}

\put(11,1){\line(1,2){2}}

\put(8,4){\line(1,2){2}}

\put(8,4){\line(1,-1){3}}

\put(10,8){\line(1,-1){3}}

\put(10,5){\line(1,-1){1}}

\put(2,8){\line(3,0){3}}

\put(8,1){\line(3,0){3}}

\put(8,4){\line(3,0){3}}

\put(10,5){\line(3,0){3}}

\put(10,8){\line(3,0){3}}

\put(5,5){\line(0,3){3}}

\put(2,5){\line(0,3){3}}

\put(8,1){\line(0,3){3}}

\put(11,1){\line(0,3){3}}

\put(10,5){\line(0,3){3}}

\put(13,5){\line(0,3){3}}

\end{picture}

\end{center}
\caption{The $Q^3$ and $G(3,2,1)$ graphs} \label{decrease}
\end{figure}
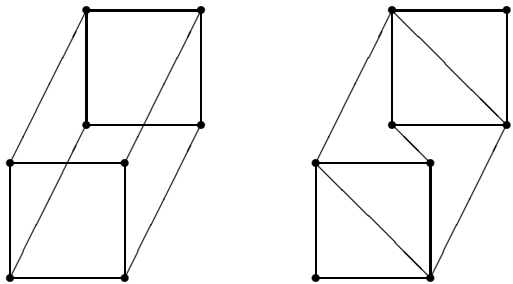

\medskip

Now in the $k$-cube there are $k\cdot 2^{k-1}$ edges and each vertex has degree $k$.  For alphabet overlap graphs, if $d=2$ and $s=1$, there are $\sim2^{k+1}$ edges, but the degree of the vertices is 2, 3, or 4.  For general values of $d$, grids have $\sim C_k\cdot k\cdot d^k$ edges (where $C_k$ is a constant depending on $k$), roughly comparable to the cardinality of the edge set of alphabet overlap graphs, but the vertex degrees of the latter graphs are around $2d$ -- smaller, in high dimensions, than the degrees $\sim2k$ of grid vertices.  {\it Alphabet overlap graphs are thus more efficient in allowing Hamiltonicity} as Proposition 1 and Theorem 2 below reveal.  The pivotal point that we wish to make is that the different architecture of alphabet overlap graphs permits Hamiltonicity in all cases, whereas this is only possible in the even case for grids.  In a similarly unexpected result, the chromatic number $\chi(G)$ of alphabet overlap graphs is shown to be larger, and often much larger, than that of bipartite grids.

The proof of Proposition 1 below is well-known (see, e.g., Theorem 5.3 in \cite{buckley} for a partial proof) but given for completeness; in addition, there are several papers, for example \cite{clark} and \cite{collins},    that focus on the nuances of the {\it rectangular} { two} dimensional case $P_m\times P_n$, with special attention being paid to the {\it numbers} of Hamiltonian cycles or paths in these graphs.
\begin{prop}
For $k\ge2$, the grid $\{1,2,\ldots,d\}^k$ permits a Hamiltonian
cycle if and only if $d$ is even, with the value of $k$ being
irrelevant.
\end{prop}

\medskip

\noindent{\bf Proof}  Assume the grid has a Hamiltonian cycle.  We must, in $d^k$ steps, travel from the ``origin" back to itself.  At each step, the $\ell^1$ (``taxicab") distance to the origin changes by $\pm1$, with the first return to the origin being at the $d^k$th step.  The numbers of $+1$ and $-1$ distance changes must thus be equal, which is possible only if $d^k$ is divisible by 2, i.e., if $d$ is even.  Conversely, let $d$ be even.  We use induction on $k$.  It is easy to verify that $\{1,2,\ldots,d\}^2$ is Hamiltonian.
Assume next that $\{1,2,\ldots,d\}^r$ is Hamiltonian and let
\begin{equation}
\{v_1,v_2,\ldots,v_{d^r},v_1\}
\end{equation}
be any Hamiltonian cycle.  Since
$\{1,2,\ldots,d\}^{r+1}$=$\{1,2,\ldots,d\}^r\times P_d$, where the $d$-path $P_d$ consists of vertices $w_1,w_2,\ldots,w_d$, we may construct a Hamiltonian cycle on $\{1,2,\ldots,d\}^{r+1}$ as follows:
\begin{eqnarray}
&&\{(v_1,w_1),(v_1,w_2),\ldots,(v_1,w_d),\nonumber\\
&&(v_2,w_d),(v_2,w_{d-1}),\ldots,(v_2,w_1),\nonumber\\
&&(v_3,w_1),(v_3,w_2),\ldots,(v_3,w_d),\nonumber\\
&&\ldots\nonumber\\
&&(v_{d^r-1},w_1),(v_{d^r-1},w_2),\ldots,(v_{d^r-1},w_d),\nonumber\\
&&(v_{d^r},w_d),(v_{d^r},w_{d-1}),\ldots,(v_{d^r},w_1),\nonumber\\
&&(v_{1},w_1)\},
\end{eqnarray}
where the $v_i$ are as in (1).  This completes the proof; note, moreover that the same argument works in any dimension for grids of the form $P_{d_1}\times\ldots\times P_{d_k}$ as long as at least one $d_i$ is even.\hfill\qed

\section{Hamiltonicity}
\begin{thm}
Alphabet overlap graphs $G(k,d,s)$ are Hamiltonian for all  $k,d\ge2$ and $s\le k-1$.
\end{thm}
%Fredricksen, Harold; Maiorana, James Necklaces of beads in $k$ colors and $k$-ary de Bruijn sequences. Discrete Math. 23 (1978), no. 3, 207--210.

\noindent We provide two proofs of this fact, the first more constructive than the second, which, furthermore works only for $s< \lfloor k/2\rfloor$ ($d$ arbitrary).  In fact, this second proof is an extension of the ``standard" proof for the case $s=1$.  Other demonstrations of this special case include an ingenious greedy algorithm proof due to Fredricksen and Maiorana \cite{hal}.  It is our hope that the elementary Proof 1 below may be of use to exhibit the existence of deBruijn cycles in other contexts.

\noindent {\bf Proof 1}  We start by proving that each alphabet overlap graph $G(k,d,1)$ is Hamiltonian.
The proof is by induction on $d$.
If $d = 1$ the graph is Hamiltonian because it has only one vertex.  Assume the claim is true when $d = n$.  Consider the alphabet overlap graph
$G(k, n+1, 1)$.  Let $L$
be an ordered set that lists the vertices along a Hamiltonian path for the
alphabet overlap graph
$G(k, n, 1)$.
Let ``$a$" be the first letter of the alphabet, and ``$b$"  the $(n+1)$st
letter of the alphabet.
Let $\Sigma$ be a string which contains {\it one and only one} $b$, that $b$ being its final
letter.  Then there are
directed edges from $\Sigma$ to the cyclic permutation of $\Sigma$ with $b$ in the
next-to-last position, and so on,
through all distinct cyclic permutations of $\Sigma$.  Now, somewhere in $L$ is the
vertex which begins
with $a$ followed by the first $k-1$ letters of $\Sigma$.  In $G(k, n+1, 1)$ there is an
edge from that vertex to
the vertex $\Sigma$.  There is also an edge between (i) the cyclic permutation
of $\Sigma$ with $b$ at
the front and (ii) the next vertex in $L$.  Therefore we can form $L'$
by inserting all of the
cyclic permutations of $s$ into $L$ as indicated.
We can do the same for each distinct string which contains one and only one
$b$, that $b$ being its
final letter.
We next employ a similar process for the insertion of  each distinct string with two $b$'s, one of which
is its final letter, not
counting strings as distinct if they are cyclic permutations of one another.
We continue in this fashion until every
string of length $k$ in an alphabet of size $n+1$ is included.  This gives a
Hamiltonian list for $G(k,
n+1, 1)$.

As an example, here is the construction of a Hamiltonian list for $G(3, 2,1)$.
We will use alphabet $\{a,b\}$.
We begin with $G(3, 1, 1)$, which has one vertex, $aaa$.  Our Hamiltonian list
has therefore a single
entry, i.e.,

\[L=\{aaa\}.\]
Now, in $G(3, 2, 1)$ there is only one vertex with only one $b$ and that in the
last position: $aab$.  The
cyclic permutations of $aab$ are: $\{aab, aba, baa,\}$, with each permutation
adjacent to the next.  We
insert these into $L$ to form $L'$:

\[L'=\{aaa,
aab,
aba,
baa\}.\]
Next, in $G(3, 2, 1)$ there are two vertices with exactly two $b$'s, one of
which is at the end.  But
since each is a cyclic permutation of the other, they are not distinct, so
we only need consider one
of them.  Choose $abb$.  The cyclic permutations of $abb$ are: $abb, bba$, and
$bab$, with each adjacent
to the next.  Following our algorithm, we insert these three strings into $L'$
between $aab$ and $aba$,
thus obtaining:

\[L''=\{aaa,
aab,
abb,
bba,
bab,
aba,
baa\}.\]

Finally, there is one vertex with three $b$'s, namely $bbb$.  Following our algorithm,
we insert this between
$abb$ and $bba$ thus obtaining the required Hamiltonian cycle

\[H=\{aaa,
aab,
abb,
bbb,
bba,
bab,
aba,
baa,aaa\}.\]

The proof of the case of general $s$ now follows easily.
The only change is that in choosing the string $\Sigma$ we only require that the
new letter be part of the
terminal string of length $s$.  Then, instead of all cyclic permutations of
$\Sigma$, we use all distinct
cyclic permutations that move the last letter of $\Sigma$ forward $k-s$ steps.  Then
we insert this list of
distinct cyclic permutations after the label with all $a$s except in the
final $k-s$ positions, those
positions agreeing with the initial segment of $\Sigma$ of length $k-s$.  The last element of this set of cyclic permutations is adjacent to $\Sigma$, and thus to the next element of the list $L$, as in the $s=1$ case.  This gives an
algorithm for
constructing a Hamiltonian path in any alphabet overlap graph $G(k, d, s)$.
\hfill\qed

\medskip

\noindent{\bf Proof 2}  We illustrate the method first for $s=1$ and $d=2$.  We wish to arrange $2^k$ binary digits in a circular array so that the set $A$ of $k$ consecutive digits has maximal cardinality $2^k$, i.e., these are all distinct.  For $k=4$ this may be accomplished as $0000101001101111$. Note that this gives the Hamiltonian cycle
$$0000\rightarrow0001\rightarrow\ldots\rightarrow1111\rightarrow1110\rightarrow1100\rightarrow1000\rightarrow0000$$
in $G(4,2,1)$.
For general values of $k$ we consider the alphabet overlap graph $G(k-1,2,1)$ as a starting point.  We draw a {\it directed} edge from vertex $(a_1,\ldots,a_{k-1})$ to vertex $(a_2,\ldots,a_{k-1},0)$, which we label $a_1\ldots a_{k-1}0$; and a directed edge from vertex $(a_1,\ldots,a_{k-1})$ to vertex $(a_2,\ldots,a_{k-1},1)$, which we label $a_1\ldots a_{k-1}1$.  Notice that we thus obtain a directed version of $G(k-1,2,1)$, with loops at the vertices $(0,0,\ldots,0)$ and $(1,1,\ldots,1)$.  Moreover, each vertex has both in- and out-degree equal to 2.  It thus (see, e.g., Theorem 1.8.1 in \cite{diestel}) has an Eulerian circuit $(e_1,e_2,\ldots, e_{2^k})$ which yields, since the edges are labeled with distinct binary sequences, the required ensemble of $2^k$ distinct sets of $k$-consecutive digits -- which we identify with a Hamiltonian cycle for $G(k,2,1)$.

When $d\ge3$, we employ the same process starting with $G(k-1,d,1)$ and observe that each vertex now has both in- and out-degree equal to $d$.  An Eulerian circuit $(e_1,e_2,\ldots, e_{d^k})$ is guaranteed to exist.  We get the required Hamiltonian cycle as before.  Finally when $\lfloor k/2\rfloor> s\ge 2$ and we are faced with showing that $G(k,d,s)$ is Hamiltonian, we consider the directed version of $G(k-s,d,s)$.  Each vertex has both in- and out-degree equal to $d^s$, so that an Eulerian circuit $(e_1,e_2,\ldots, e_{d^k})$ exists, where an edge between vertices $(a_1,\ldots,a_{k-s})$ and $(a_{s+1},\ldots,a_{k})$ is denoted by $(a_1\ldots a_{k})$.  The rest of the argument is as before.\hfill\qed

\section{Chromatic Number of Alphabet Overlap Graphs}
In this section and the next we will often use symbol $t$ to denote the ``tag length" $k-s$.  This is particularly convenient when $s\ge k/2$, when we can think of a word of length $k$ as being of the form $t_1vt_2$, where (i) $t_1,t_2$ are tags of length $t$, and (ii) $v$, possibly the empty word, has length $k-2t$.

\begin{thm}If  $t\le k/2$, i.e. if $s\ge k/2$, then the chromatic number of $G(k,d,s)$ is given by $\chi(G(k,d,s))= d^{k-2t} + d^t$.
 \end{thm}

\medskip

\noindent {\bf Proof} Let $w^*$ be any fixed word of length $k-2t$
(if $k-2t$ is zero, then $w^*$ is the empty word).  $G(k,d,s)$ has
an induced subgraph consisting of all words of the form
$t_1w^*t_2$, where the values of $t_1$ and $t_2$ range over all
tags with $t_1\ne t_2$.  If we order the $d^t$ tags,
lexicographically for example, we can arrange these words in a
square array, in the form of a matrix missing its main diagonal,
so that the word $t_iw^*t_j$ is in the $i$th row and $j$th column.
Given any word in this array, we can find all adjacent words in
the array as follows.  The word symmetric about the main diagonal
is adjacent, as are all words on the same row or the same column
with the symmetric word.  For the purposes of this paper, call any
graph whose vertices are $v_{ij}, i \ne j, 1\le i,j\le n, n\ge2$,
with an edge between $v_{ij}$  and $v_{xy}$ iff $x = j$ or $i =
y$, an {\it alphabet overlap (AO-) matrix graph of order $n$}. For
any $n$ we obtain an AO matrix graph of order $n$ from $G(2,n,1)$. 
The chromatic number of any AO matrix graph of order $n$ is 
$n$, as we prove in the following paragraphs.

Since no two words in the same column are adjacent, we can color an AO matrix graph of degree $n$ with $n$ colors, by coloring every word in a given column the same color, while coloring words in different columns different colors.  We use induction to prove that this is a minimal coloring.

The AO matrix graph of order 3 has chromatic number at least 3
since it contains a triangle, and chromatic number at most 3
because we can color all entries in each column the same color. By
inspection, the only minimal colorings of this graph either have
all entries in each column the same color or else have all the
entries in each row the same color.  Also by inspection, the AO
matrix graph of order 4 has chromatic number exactly 4, and the
only minimal colorings have either monochromatic columns or
monochromatic rows. Assume that the AO matrix graph of order $n =
N \ge4$ has chromatic number $N$ and that the only minimal colorings
have either monochromatic columns or monochromatic rows. Consider
a minimal coloring of the AO matrix graph of order $N+1$. We know
that $N+1$ colors suffice. Is an $N$ coloring possible? The upper
left hand $N\times N$ corner satisfies the induction hypothesis,
and so requires $N$ colors. Further, for any $N$-coloring of the
upper left hand corner, either all entries in each row are the
same color or all entries in each column are the same color.
Assume, without loss of generality, that all entries in each
column are the same color. The same is true of the lower right
hand $N\times N$ corner, which in the case of $N\ge 4$ overlaps
the upper left hand corner in at least two places in each of the
middle $N-1$ columns, and so must have monochromatic columns
rather than monochromatic rows. Therefore the middle $N-1$ columns
are monochromatic, as are the first column, excepting possibly its
last entry, and the last column, excepting possibly its first
entry.

There are four sets of colors to be considered, namely the color of the entries in the last column (except for the first entry); the color of the entries in the first column (except for the last entry); the color of the entry in the $(1,N+1)$ position; and the color of the entry in the $(N+1,1)$ position.  Call these $\chi_1, \chi_2, \chi_3$ and $\chi_4$ respectively.  Now {if} $\chi_1=\chi_2$, it is easy to verify that $\chi_3$ and $\chi_4$ must be two {new colors} forcing us to use $N+2$ colors.  We thus set $\chi_1\ne\chi_2$.  From this it follows that $\chi_1=\chi_3\ne \chi_2=\chi_4$ for a total of $N+1$ colors  with each column monochromatic.  This establishes that the chromatic number of AO matrix graphs is equal to the number of columns.

Now fix a tag $t^*$, and consider the clique of all words of the form $t^* m_it^*$.  All of the $d^{k-2t}$ words of this form are adjacent to one word in every column of each AO matrix graph that is an induced subgraph of $G(k,d,s)$.  Thus none of them can be colored any of the $d^t$  colors used in an AO matrix graph.  Since they form  a clique, we must color them all different colors.  Therefore $G(k,d,s)$ requires at least $d^{k-2t} + d^t$ colors.

To see that this number of colors suffices, given any word $t_am_bt_c$ , if $t_a=t_c\ne t^*$, color that word the same color as $t^*m_bt^*$.  Otherwise, color it the same color as $t_aw^*t_c$.  If $w^*$ is the empty word, we have already colored $t_am_bt_c=t_at_c$ in the AO matrix graph.  \hfill\qed

\begin{thm}
If $t > k/2$, then $\chi(G(k,d,s)) \le 1 + d^t=1+d^{k-s}$.
\end{thm}

\medskip

\noindent{\bf Proof}
Let $t > k/2$. There is an isomorphism between $G(k,d,s)$ and a subgraph of $G(2t, d, t)$ under which vertex $xmz$ corresponds to vertex $xmmz$, where $x$ and $z$ are words of length $k - t=s$ and $m$ is a word of length $2t - k$, and where $xm$ and $mz$ are tags both in $G(k,d,s)$ and in $G(2t,d,t)$.   In other words, the induced subgraph of $G(2t,d,t)$ isomorphic to $G(k, d,s)$ is the graph whose vertices have the form $xmmz$.  We call a graph having this form a {\it reduced AO matrix graph.}  Since in $G(2t,d,t)$ each tag is exactly half of each word, Theorem 3 applies and we conclude that the chromatic number of $G(2t,d,t)$ is $d^{2t - 2t} +d^t = 1 + d^t$, which establishes the theorem.\hfill\qed

\medskip

\noindent{\bf Remark}
We can improve on the upper bound in Theorem 4 using the following algorithm:

Consider the columns of the reduced AO matrix graph for $G(2t, d,
t)$.  We can color two of these columns the same color if they
contain no adjacent pairs of entries.  The column containing words
ending in $t_i$ and the column containing words ending in $t_j$
have adjacent entries if and only if either the word $t_it_j$ or
the word $t_jt_i$ is in the subgraph.    This happens if and only
if the last $2t - k$ letters in $t_i$ match the first $2t - k$
letters in $t_j$ or vice versa.  This happens if and only if $t_i$
and $t_j$ are adjacent in $G(t,d,k-t)$. Thus there is a one-to-one
correspondence between the colorings of the columns of the reduced
AO matrix graph of $G(2t,d,t)$ and the colorings of $G(t,d,k-t)$.
If $k-t\ge t/2$, that is if $t\le\frac{2}{3}k$, then by Theorem 3,
the chromatic number of $G(t,d,k-t)$ is $d^{2k-3t}+d^{2t-k}$, so
the columns of the reduced AO matrix graph of $G(2t,d,t)$ can be
colored with that number of colors. Therefore $G(k,d,s)$ can be
colored with $1+d^{2k-3t}+d^{2t-k} $ colors, a better bound on the
chromatic number than the one in Theorem 4. For example, Theorem 4
gives $\chi(G(5,2,2))\le 9$, while the algorithm outlined in this
remark reduces the upper bound to $\chi(G(5,2,2))\le 5$.

If $t >\frac{2}{3} k$, we can repeat the process.  Eventually we find an upper bound on the chromatic number better than the upper bound in Theorem 4.  We do not get an exact value for the chromatic number by this method, however, because reduced AO matrix graphs, unlike AO matrix graphs, may have minimal colorings in which neither the columns nor the rows are monochromatic.

\section{Planarity and domination number of AO graphs}
In this section, we provide some baseline results on planarity and   domination numbers of alphabet overlap graphs.
\begin{thm}
If $t\le k/2$ the only non-trivial planar AO-graphs are when $d = 2, 3$, $t = 1$, and $k=2$.
\end{thm}

\medskip

\noindent{\bf Proof} If we have 4 distinct tags $\alpha,\beta,\gamma,\delta$, then we can construct a $K_{3,3}$ subgraph as follows.  Let $w$ be any word of length $k - 2t$ (possibly the empty word).  Then $\alpha w\beta, \alpha w\gamma,\alpha w\delta, \beta w\alpha, \gamma w\alpha, \delta w\alpha$ form a non-planar $K_{3,3}$ subgraph. If there are not 4 distinct tags, then $d = 2$ or 3 and $t$ = 1.  Now if $d=2$ it is east to verify that $G(2,2,1)$ {\it is} planar but that $G(3,2,2)$ {\it isn't} since it contains a bipartite $K_{4,4}$ subgraph with classes $\{101,111,010,000\}$ and $\{110,100,001,011\}$ (there are additional edges between vertices in the same color classes).  In a similar fashion, it is not too hard to draw a planar version of $G(2,3,1)$.  Since for $a<b$, $G(k,a,s)$ is a subgraph of $G(k,b,s)$, this completes the proof.\hfill\qed

\begin{thm}
If $t\le k/2$, then the domination number of $G(k,d,s)$ is $d^t$.
\end{thm}

\medskip\noindent{
\bf Proof} Clearly, if $x$ is any word of length $k - t=s$, and if we name the $n = d^t$ tags  $t_1,\ldots,t_n$, then the words $t_ix$ form a dominating set.  Now, suppose we have a minimal dominating set $S$ in which some tag, say $t_1$, does not appear at the beginning of any word in $S$.  Consider the set of all words of the form $t_iwt_1$, for any fixed word $w$ of length $k - 2t$.  For every tag $t_i$ there must be some word in $S$ that ends in $t_i$, so there are $d^t$ vertices in $S$.\hfill\qed

\section{Open Problems} We mention two open problems.  First, it would be most interesting, in line with the investigations in \cite{clark} and \cite{collins}, to estimate the {\it number} of Hamiltonian cycles in alphabet overlap graphs $G(k,d,s)$; this number is known when $s=1$. Second, we feel that several structural properties of $G(k,d,s)$ are worth studying.  These might include
\begin{enumerate}
\item connectivity properties;
\item distance properties;
\item existence of cycles of various lengths (e.g., Is $G(k,d,s)$ pancyclic?);
\item colorings and cliques (e.g., What is $\chi(G(k,d,s)$ when $s=1$?  What are the exact values of $\chi$ in general?  What is the clique number of $G(k,d,s)$?); and
\item the existence of special substructures.
\end{enumerate}
\medskip

\noindent{\bf Acknowledgment}  The research of the first named author was  supported by NSF Grants DMS-0049015 and DMS-0139291.  The research of both the first and second named authors was supported by NSF Grant DMS-0122278.

\end{document}